\documentclass{article}
\usepackage{colt10e}

\usepackage{natbib,amsmath,amssymb}

\usepackage{bm,amsfonts}

\usepackage{algorithm,algorithmic}

\DeclareMathOperator*{\argmax}{\rm argmax}
\newcommand{\eqdef}{\equiv}
\newcommand{\so}{\mathrm{o}}
\newcommand{\lo}{\mathrm{O}}
\newcommand{\ep}{\epsilon}
\newcommand{\de}{\delta}
\newcommand{\id}{ \mbox{\rm 1}\hspace{-0.63em}\mbox{\rm \small 1\,}}
\newcommand{\idx}[1]{\id\left[#1\right]}
\newcommand{\idxx}[1]{\id\!\left[#1\right]}
\newcommand{\E}{\mathrm{E}}
\newcommand{\dmin}{D_{\mathrm{inf}}}
\newcommand{\hen}[2]{\frac{\partial #1}{\partial #2}}
\newcommand{\n}{\nonumber}
\newcommand{\nn}{\nonumber\\}
\newcommand{\scap}{\,\cap\,}
\newcommand{\fhat}{\hat{F}}
\newcommand{\muhat}{\hat{\mu}}
\newcommand{\fhatn}[2]{\fhat_{#1}(#2)}
\newcommand{\fhatt}[2]{\fhat_{#1,#2}}
\newcommand{\muhatn}[2]{\muhat_{#1}(#2)}
\newcommand{\muhatt}[2]{\muhat_{#1,#2}}
\newcommand{\rd}{\mathrm{d}}

\newtheorem{remark}[theorem]{Remark}
\newtheorem{proposition}[theorem]{Proposition}
\newenvironment{proof1}{\noindent{\bf Proof:}}{\vspace{1.5mm}}
\newenvironment{proof2}[1]{\noindent{\bf Proof #1:}}{\vspace{1.5mm}}
\newcommand{\fqed}{\hfill $\blacksquare$ \par}
\def\dqed{\relax\tag*{\qed}}

\newcommand{\mom}{\frac{1}{1-\mu}}
\newcommand{\moms}{\mbox{$\mom$}}
\newcommand{\momss}{(1-\mu)^{-1}}

\newcommand{\supp}{\mathrm{supp}}
\newcommand{\di}[1]{\dmin(F_{#1},\mu^*)}
\newcommand{\jp}[2]{J'_{#1}(#2)}
\newcommand{\logn}{\frac{\log n}{(1-\ep)(1-\para)\di{i}}}
\newcommand{\nus}{\nu^*}
\newcommand{\la}{\lambda}
\newcommand{\La}{\mathrm{\Lambda}}
\newcommand{\e}{\mathrm{e}}

\newcommand{\com}{\enspace ,}
\newcommand{\per}{\enspace .}
\newcommand{\uI}{\tilde{\La}^*} 

\newcommand{\cala}{\mathcal{A}}
\newcommand{\iopt}{\mathcal{I}_{\mathrm{opt}}}
\renewcommand{\j}[2]{J_{#1}{(#2)}}

\newcommand{\levy}{d_{\mathrm{L}}}
\newcommand{\bbR}{\mathbb{R}}
\newcommand{\calX}{\mathcal{X}}
\newcommand{\calV}{\mathcal{V}}
\newcommand{\Cb}{C_{\mathrm{b}}}
\newcommand{\calE}{\mathcal{E}}
\newcommand{\para}{r}
\newcommand{\tmu}{\mu'+\de}
\newcommand{\dmint}{L_{\max}} 
\newcommand{\cc}{c_0}
\newcommand{\cd}{\frac{1-\E(F)}{1-\mu}}
\newcommand{\since}[1]{\quad\left(\mbox{#1}\right)}
\newcommand{\mg}{\La^*}
\newcommand{\mut}{\mu}
\newcommand{\mt}{\tilde{\mu}}
\newcommand{\expn}[1]{\e^{-#1}}
\newcommand{\expp}[1]{\e^{#1}}

\title{
Finite-time Regret Bound of a Bandit Algorithm for the Semi-bounded Support Model
}

\author{Junya Honda and Akimichi Takemura\\
The University of Tokyo, Japan.
\\\texttt{\small \{honda,takemura\}@stat.t.u-tokyo.ac.jp}
}

\begin{document}

\maketitle
\allowdisplaybreaks[1]
\begin{abstract}
In this paper we consider stochastic multiarmed bandit problems.
%
Recently a policy, DMED, is proposed and proved to achieve the asymptotic
bound for the model that each reward distribution is supported in a known
bounded interval, e.g.~$[0,1]$.
However, the derived regret bound is described in an asymptotic form and
the performance in finite time has been unknown.
We inspect this policy and derive a finite-time regret bound
by refining large deviation probabilities to a simple finite form.
Further, this observation reveals that
the assumption on the lower-boundedness of the support is not essential
and can be replaced with a weaker one,
the existence of the moment generating function.

\end{abstract}

\section{Introduction}\label{section-intro}
In the multiarmed bandit problem a gambler pulls arms of a slot machine
sequentially so that the total reward is maximized.
There is a tradeoff between exploration and exploitation
since he cannot know the most profitable arm
unless pulling all arms infinitely many times.

There are two main formulations for this problem:
stochastic and nonstochastic bandits.
In the stochastic setting rewards of each arm follow an unknown distribution
\citep{gittins,conti1,vermorel},
whereas the rewards are detemined by an adversary
in the nonstochastic setting \citep{adversarial}.
In this paper we consider the stochastic bandit, where rewards
of arm $i\in\{1,\dots,K\}$
are i.i.d.~sequence from unknown distribution $F_i\in\mathcal{F}$ with expectation $\mu_i$
for a model $\mathcal{F}$ known to the gambler.
For the maximum expectation $\mu^*\equiv \max_i \mu_i$,
we call an arm $i$ optimal if $\mu_i=\mu^*$ and suboptimal otherwise.
If the gambler knows each $\mu_i$ beforehand, it is best to
choose optimal arms 
at every round.
A {\it policy} is a strategy of the gambler for choosing arms based on
the past result of plays.
The performance of a policy is measured by the loss called
expected regret or {\it regret}, in short, 
given by
\begin{eqnarray}
\sum_{i:\mu_i<\mu^*}(\mu^*-\mu_i)\E[T_i(n)]\com\n
\end{eqnarray}
where $T_i(n)$ is the number of plays of arm $i$ through the first $n$ rounds.
Since we regard each $\mu_i$ as a unknown constant fixed in advance,
we consider how we can reduce $\E[T_i(n)]$ for each suboptimal arm $i$
to achieve a small regret.

\citet{robbins} first considered this setting
and \citet{lai} gave a framework for determining an optimal policy
by establishing a theoretical bound for the regret.
Later this theoretical bound was extended to multiparameter or nonparametric models
$\mathcal{F}$ by \cite{burnetas}.
In their paper, it was proved that any policy satisfying a mild regularity condition
satisfies
\begin{eqnarray}
\E[T_i(n)]\le \frac{1-\so(1)}{\dmin(F_i,\mu^*;\mathcal{F})}\log n\com \label{bound_burnetas}
\end{eqnarray}
where $\dmin(F,\mu;\mathcal{F})$ is defined in terms of
Kullback-Leibler divergence $D(\cdot\Vert\cdot)$ by
\begin{eqnarray}
\dmin(F_i,\mu^*;\mathcal{F})=\inf_{G\in\mathcal{F}:\E_{G}[X]> \mu^*}D(F_i\Vert G)\per\n
\end{eqnarray}

The most popular model in the nonparametric setting is the family
of distributions with supports contained
in a known bounded interval, e.g.~$[0,1]$. 
For this model, which we denote by $\cala_0$,
it is known that fine performance can be obtained by
policies called Upper Confidence Bound (UCB) \citep{ucb,ucbv,kl_ucb}.
However, although some bounds for regrets of UCB policies
have been obtained in a non-asymptotic form,
they do not necessarily achieve the asymptotic theoretical bound.
%

Recently \citet{honda_colt} proposed
Deterministic Minimum Empirical Divergence (DMED) policy,
which chooses arms based on an index $\dmin(\fhat_i,\mu;\cala_0)$,
or simply written as $\dmin(\fhat_i,\mu)$,
for empirical distribution $\fhat_i$ of arm $i$.
Whereas DMED achieves
the theoretical bound asymptotically, the evaluation
heavily depends on an asymptotic analysis
and any finite-time regret bound has been unknown.
Further, in the analysis of DMED, the assumption on the lower bound of the support
seems to be a technical one needed for the proof.
For example, 
the gambler does not have to know that the lower bound of the support is zero
if he knows that the upper bound is one.
%

{\bf Our Contribution.}
Based on the above observation,
we consider the family $\cala$ of distributions on $(-\infty,1]$ instead of $\cala_0$.
We first show that $\dmin(F,\mu;\cala_0)=\dmin(F,\mu;\cala)$ for all $F\in \cala_0$.
Thus,
although the gambler has more candidates for the true distribution
of each arm in the model $\cala$ than in $\cala_0$,
the theoretical bound \eqref{bound_burnetas} does not vary
between $\cala_0$ and $\cala$.

Next we provide a finite-time regret bound of DMED for all distributions in $\cala$
with moment generating functions existing in some neighborhood of the origin.
Since nonstochastic bandits inevitably require the boundedness of the support,
we can now assert that an advantage of assuming stochastic bandits is that the
semi-bounded rewards can be dealt with in the nonparametric setting.

{\bf Technical Approach.}
In the evaluation of DMED
it is essential to evaluate the probability that
$\dmin(\fhat_i,\mu)$ deviates from $\dmin(F_i,\mu)$.
Note that for policies based on the index $\dmin(\fhat_i,\mu)$,
finite-time regret bounds have been
derived for the case that each distribution is supported in a {\it finite} subset of $[0,1]$
\citep{maillard,honda_ML}.
The advantage of assuming finiteness is that
Sanov's theorem gives a non-asymptotic large deviation probability. 
However the regret bounds derived by this technique contain
a finite but exceedingly large term
\begin{eqnarray}
\sum_{t=1}^{\infty}t^{|\supp(F_i)|}\,\expn{at}\com\n
\end{eqnarray}
where $|\supp(F_i)|$ denotes the size of the support of  $F_i$ and
the polynomial $t^{|\supp(F)|}$ appears as a total number
of possible empirical distributions from $t$ samples from $F_i$.
Similarly, whereas non-asymptotic Sanov's theorem is also known
for continuous support distributions (see \citet[Ex.\,6.2.19]{LDP}),
it requires the total number of $\ep$-balls to cover a set of distributions
as a coefficient.
Thus,
although it is not impossible to derive a finite-time regret bound
by a naive application of the non-asymptotic Sanov's theorem,
it becomes very complicated and unrealistic.

To avoid counting or covering the possible empirical distributions,
we exploit the following fact
\begin{eqnarray}
\dmin(\fhat_i,\mu)=\max_{0\le\nu\le\mom}\E_{\fhat_i}[\log(1-(X-\mu)\nu)]\per
\label{intro_dmint}
\end{eqnarray}
Although it involves a maximization operation,
it is merely an empirical mean of random variables $\log(1-(X_t-\mu)\nu)$
where each $X_t$ follows distribution $F_i$.
By Cram\'er's theorem we can bound the large deviation probability
for such a finite dimensional empirical mean
by an exponential function with a simple coefficient.

Another difficulty for our setting is that $\dmin(F,\mu)=\dmin(F,\mu;\cala)$ is
neither bounded nor continuous in $F\in\cala$ unlike the case of $\cala_0$,
which makes the evaluation of the exponential rate
for the large deviation probability of $\dmin(\fhat_i,\mu)$ much harder.
The key to this problem also lies in \eqref{intro_dmint}.
Since it is an expectation of a logarithmic function on $X$,
the effect of the tail weight is weaker than the polynomial function $X^1=X$.
Thus the large deviation probability of the joint distribution of
$(\dmin(\fhat_i,\mu),\E_{\fhat_i}[X])$
can be evaluated on the same regularity condition
as that for the empirical mean $\E_{\fhat_i}[X]$ alone,
namely, the existence of the moment generating function of $F_i$
in some neighborhood of the origin.

{\bf Paper Outline.}
In Sect.\,\ref{section-preliminary} we give definitions used
throughout this paper
and introduce DMED policy proposed for distributions on $[0,1]$.
In Sect.\,\ref{section_main}, we give the main results of this paper
on the finite-time regret bound of DMED for distributions on $(-\infty,1]$.
The remaining sections are devoted to the proof of the main results.
We extend some results for the support $[0,1]$ to $(-\infty,1]$
in Sect.\,\ref{section_extension}.
We derive a large deviation probability for $\dmin(F,\mu)$
in a non-asymptotic form in Sect.\,\ref{section_LDP}.
We conclude this paper 
in Sect.\,\ref{section_remarks}.
We give some results on large deviation principle in Appendix \ref{append_ldp}.
A proof of the main theorem is given in Appendix \ref{proof_opt}.

\section{Preliminaries}\label{section-preliminary}
Let $\cala_a,\,a\in(-\infty,1),$ be the family of probability distributions on $[a,1]$.
We denote the family of distributions on $(-\infty,1]$ by
$\cala_{-\infty}$ or simply $\cala$.
For $x\in\bbR$ and $F\in \cala$,
the cumulative distribution is denoted by $F(x)=F((-\infty,x])$.
For the metric of $\cala_a$ we use L\'evy distance
\begin{eqnarray}
\levy(F,G)\equiv \inf\{h>0: F(x-h)-h\le G(x)\le F(x+h)+h\}\per\n
\end{eqnarray} 
$\E_F[\cdot]$ denotes the expectation under $F\in\cala$.
When we write e.g.\,%
$\E_F[u(X)]$ for a function $u:\mathbb{R}\to\mathbb{R}$,
$X$ denotes a random variable with distribution $F$.
The expectation of $F$ is denoted by $\E(F)\equiv \E_F[X]$.
We always assume that
the moment generating function $\E_{F}[\expp{\la X}]$ is finite
in some neighborhood of the origin $\la=0$.

Let $T_i(n)$ be the number of times that arm $i$ has been pulled through the first $n$ rounds.
$\fhatt{i}{t}$ and $\muhatt{i}{t}$ denote the empirical distribution and the mean of arm $i$
when arm $i$ is pulled $t$ times.
$\fhatn{i}{n}\equiv \fhatt{i}{T_i(n)}$ and $\muhatn{i}{n}\equiv \muhatt{i}{T_i(n)}$
denote the empirical distribution and the mean of arm $i$ at the $n$-th round.
The largest empirical mean after the first $n$ rounds is denoted
by $\hat{\mu}^*(n)\equiv \max_i\muhatn{i}{n}$.



In this paper we analyze DMED policy proposed by \cite{honda_colt}.
It is described as Algorithm \ref{alg_dmed}, where
\begin{eqnarray}
\dmin(F,\mu; \cala_a)\equiv \inf_{G\in \cala_a: \E(G)> \mu}D(F \Vert G)\per
\label{def_dmin}
\end{eqnarray}
\begin{algorithm}[t]
\caption{DMED Policy}\label{alg_dmed}
\begin{algorithmic}
\sonomama{Parameter: $\para\in(0,1).$}
\sonomama{Initialization: $L_C,L_R:=\{1,\cdots,K\},\,L_N:=\emptyset,\,n:=K$.
\textnormal{Pull each arm once}.}
\sonomama{Loop:}
\STATE
\begin{description}
\setlength{\parskip}{-1cm}
\item[1.] For $i\in L_C$ in ascending order,
\begin{description}\setlength{\parskip}{0cm}
\item[1.1.] $n:=n+1$ and pull arm $i$. $L_R:=L_R\setminus\{i\}$.
\item[1.2.] $L_N:=L_N\cup \{j\}$ 
for all
$j\notin L_R$
 such that the following $\jp{n}{j}$ occurs:
 \begin{eqnarray}
\jp{n}{j}\,\equiv\,
\{(1-r)T_i(n)\dmin(\fhatn{i}{n},\muhat^*(n);\cala_a)
\le \log n\}.\label{def-jn}
\end{eqnarray}
\end{description}
 \item[2.] $L_C,L_R:=L_N$ and $L_N:=\emptyset$.
\end{description}
\end{algorithmic}
\end{algorithm}%
%
Note that this policy is parametrized by $r\in(0,1)$ in this paper,
which was fixed to $r=0$ in the original proposal.
This parameter arises because
some properties on $\dmin(F,\mu;\cala_a)$, such as boundedness and continuity,
do not hold for $a=-\infty$.
For $r>0$ we conservatively (i.e.~more often) choose seemingly suboptimal arms.
As a result, the coefficient of the logarithmic term becomes
$1/(1-r)$ times the theoretical bound.

Another minor change is that
$\log n$ in \eqref{def-jn} was $\log n -\log T_i(n)$ in the original proposal.
It is described in \citet{honda_colt}
that the term $\log T_i(n)$ is only for improvement of simulation results
and has no importance for the asymptotic analysis.
In this paper we avoid this term since it makes the constant term in the finite-time
analysis much more complicated.

For the setting of $a=0$,
the regret of DMED is evaluated as follows.
\begin{proposition}[{\citet[Theorem 4]{honda_colt}}]\label{dmin_colt}
Let $\ep>0$ be arbitrary.
Under DMED policy with $r=0$,
it holds for all $(F_1,\dots,F_K)\in\cala_0^K$ and suboptimal arms $i$ that
\begin{eqnarray}
\E
[T_i(n)]\le
\frac{1+\ep}{\dmin(F_i,\mu^*; \cala_0)}\log n +\lo(1).\n
\end{eqnarray}
\end{proposition}
This bound is asymptotically optimal in view of the theoretical bound
\eqref{bound_burnetas}.

Now define
\begin{eqnarray}
L(\nu; F,\mu)&\eqdef&\E_F[\log(1-(X-\mu)\nu)]\com\nn
\dmint(F,\mu)&\equiv&\max_{0\le\nu\le \mom}L(\nu;F,\mu)\per\label{dmint}
\end{eqnarray}
Functions $L$ and $\dmint$ correspond to the Lagrangian function
and the dual problem of $\dmin(F,\mu; \cala_a)$, respectively. 
\begin{proposition}[{\citet[Theorem 5]{honda_colt}}]\label{prop_dmin}
For all $F\in\cala_0$ and $\mu<1$ it holds that
$\dmin(F,\mu;\cala_0)=\dmint(F,\mu)$.
\end{proposition}

\section{Main Results}\label{section_main}
We now state the main result of this paper in Theorems \ref{thm_bound} and \ref{optMD}.
We show that the theoretical bound does not depend on knowledge of the
lower bound of the support in Theorem \ref{thm_bound}
and that the theoretical bound is actually achievable by DMED in Theorem \ref{optMD}.
\begin{theorem}\label{thm_bound}
Let $a\in [-\infty,1)$ and $F\in \cala_a$ be arbitrary.
{\rm (i)}
$\dmin(F,\mu;\cala_a)=\dmin(F,\mu;\cala)$. 
{\rm (ii)} If $\mu<1$ then $\dmin(F,\mu;\cala)=\dmint(F,\mu)$.
\end{theorem}
We prove this theorem in the next section.
The part (i) of this theorem means that
the theoretical bound does not depend on whether
we know that the support of distributions is bounded from below by $a$
or
we have to consider the possibility that the support of distributions may not
be lower-bounded.
Furthermore, from (ii), we can express the theoretical bound
in the same expression as $\cala_0$
for any distribution in $\cala$.
In view of this theorem we sometimes write $\dmin(F,\mu)$
instead of more precise $\dmin(F,\mu;\cala_a)$
or $\dmin(F,\mu;\cala)$.

Let $\iopt\equiv \{i:\mu_i=\mu^*\}\subset\{1,\cdots,K\}$ be the set of optimal arms
and $\mu'\equiv\max_{i\notin\iopt}\mu_i$ be the second optimal expected value.
Define Fenchel-Legendre transform of the moment generating function
of $F_k$ as
\begin{align}
\mg_k(x)\equiv\sup_{\la\in \bbR}\{\la x-\log\E_{F_k}[\expp{\la X}]\}\per\label{def_ik}
\end{align}
Then $\E[T_i(n)]$ is bounded
for $\xi_{i,\ep,\de}\equiv \ep \dmin(F_i,\mu^*)-\de/(1-\mu^*)$
as follows.
\begin{theorem}\label{optMD}
Assume that $\mu^*<1$.
Let $\ep>0$ and $i\notin\iopt$ be arbitrary
and fix any $\de\in(0,\mu^*-\mu')$ such that
$\xi_{i,\ep,\de}>0$.
Then for all $n>0$ 
\begin{align}
\E[T_i(n)]\le \logn+ C\com\n
\end{align}
where,
for $\uI(\cdot,\cdot,\cdot)$ defined in \eqref{ui_def},
the constant term is given by
\begin{align}
C&=
\frac{1}{1-\expn{\uI(\xi_{i,\ep,\de}\,,\,\mu_i,\,\mu^*)}}
+\!\sum_{k\in\iopt} \frac{K}{1-\expn{\mg_k(\mu^*-\de)}}
+\!\sum_{k\notin\iopt} \frac{K}{1-\expn{\mg_k(\tmu)}}\nn
&\quad
+\min_{k\in\iopt}\!\!\left\{
\frac{2(1+K)}{1-\expn{\mg_k(\tmu)}}
+\frac{2\e}{r\big(1-\expn{r \mg_k(\tmu)}\big)^2}
\right\}
\per\n
\end{align}
\end{theorem}
We prove this theorem in Appendix \ref{proof_opt}.
The proof is largely the same as that of \citet[Theorem 4]{honda_colt},
with difference that asymptotic large deviation probabilities are replaced with
non-asymptotic forms in Theorems \ref{ldp_upper} and \ref{ldp_lower}.

As described in Prop.\,\ref{cramer} (iii) of Appendix \ref{append_ldp},
$\mg_k(\cdot)$ corresponds to the exponential rate of the probability
on the sample size that the empirical mean of arm $k$ deviates from
its expectation.
We can bound this rate in an explicit form for some cases.
For example, it can be bounded by the variance for the case that
the support of $F_k$ is bounded from below \citep[Theorem 1]{hoeffding}.
However, it seems to be impossible to bound
the rate by its finite-degree moments for an optimal arms $k\in\iopt$ in general case,
although it is possible for suboptimal arms $k\notin \iopt$
\citep[Theorem 3]{hoeffding}.

\begin{remark}
The derived bound is somewhat weaker than that for the bounded support model
in Prop.\,\ref{dmin_colt} since the bound in this theorem
contains the coefficient $1/(1-r)$ in the logarithmic term.
We can remove the effect of the parameter $r$ from the logarithmic term
by letting $r$ depend on $T_i(n)$, e.g., $r=1/\sqrt{T_i(n)}$.
However, it makes the analysis longer and we omit
the evaluation of this version for lack of space.
\end{remark}

\section{Properties of $\dmin$ in the Semi-bounded Support Model}\label{section_extension}
In the analysis of DMED it is essential to investigate the function $\dmin(F,\mu;\cala)$.
In this section we extend some results on $\dmin(F,\mu;\cala_0)$ in \cite{honda_colt} for
our model $\cala=\cala_{-\infty}$ and prove Theorem \ref{thm_bound}. 
 
First we consider the function $L(\nu;F,\mu)=\E_F[\log (1-(X-\mu))\nu]$.
The integrand $l(x,\nu)\equiv\log (1-(x-\mu)\nu)$ is differentiable
in $\nu\in (0,\momss)$ for all $x\in (-\infty,1]$ with
\begin{eqnarray}
\frac{\partial l(x,\nu)}{\partial\nu}=-\frac{x-\mu}{1-(x-\mu)\nu}\com\qquad
\frac{\partial^2 l(x,\nu)}{\partial \nu^2}=-\frac{(x-\mu)^2}{(1-(x-\mu)\nu)^2}\per\n
\end{eqnarray}
Since they are bounded in $x\in (-\infty, 1]$, the integral $L(\nu;F,\mu)$
is differentiable in $\nu$ with
\begin{eqnarray}
 L'(\nu;F,\mu)&\eqdef&\hen{L(\nu;F,\mu)}{\nu}
=-\E_F\left[\frac{X-\mu}{1-(X-\mu)\nu}\right]\com\nn
 L''(\nu;F,\mu)&\eqdef&\hen{^2L(\nu;F,\mu)}{\nu^2}
=-\E_F\left[\frac{(X-\mu)^2}{(1-(X-\mu)\nu)^2}\right]\per\n
\end{eqnarray}
From these derivatives the optimal solution
$\nu^*(F,\mu)\equiv\argmax_{0\le\nu\le\momss} L(\nu;F,\mu)$
of \eqref{dmint} exists uniquely and satisfies the following lemma.
\begin{lemma}\label{lem_integrable}
Assume that $\E(F)\le \mu<1$ holds.
If $\E_F[(1-\mu)/(1-X)]\le 1$ then $\nu^*(F,\mu)=\break \momss$ and
therefore $\E_F[1/(1-(X-\mu)\nu^*)]\le 1$.
Otherwise, $L'(\nu^*;F,\mu)=0$ and \break
$\E_F\left[1/(1-(X-\mu)\nu^*)\right]=1$.
\end{lemma}
The differentiability of $\dmint(F,\mu)$ in $\mu$ also holds as in the case of bounded support.
\begin{lemma}\label{lem_bibun}
For $\mu>\E(F)$, $\dmin(F,\mu)$ is differentiable with
\begin{eqnarray}
\frac{\rd \dmin(F,\mu)}{\rd \mu}=\nu^*(F,\mu)\le \mom\per\n
\end{eqnarray}
\end{lemma}
We omit the proofs of Lemmas \ref{lem_integrable} and \ref{lem_bibun}
since they are the same as Theorems 3 and 5 of \citet{honda_ML} where
the assumption on the support is not exploited.

Define $F_{(a)}\in \cala_a$
as the distribution obtained by transferring the probability of $(-\infty,a)$
under $F$ to $x=a$, that is,
\begin{eqnarray}
F_{(a)}(x)\equiv\begin{cases}
0 &x<a\com\\
F(x)&x\ge a\per
\end{cases}\n
\end{eqnarray}
Now we give the key to extension for the semi-bounded support in the following lemma,
which shows that the effect of the tail weight is bounded uniformly if the expectation
is bounded from below.
\begin{lemma}\label{lem_fa}
Fix arbitrary $\mu,\mt<1$ and $\ep>0$. 
Then there exists $a(\ep)$ such that
$|\dmint(F_{(a)},\mu)\allowbreak -\dmint(F,\mu)|\le \ep$ for all $a\le a(\ep)$ and
$F\in\cala$ such that $\E(F)\ge \mt$ .
\end{lemma}
\begin{proof1}
Take sufficiently small $a<\min\{0,\mu\}$
and define $A=(-\infty,a), B=[a,1]$.
Note that $F(A)+F(B)=1$.
First we have
\begin{eqnarray}
F(A)&\le&\frac{1-\mt}{1-a}\label{bound_ga}\\
\int_A x\rd F(x)&\ge&\mt-1+F(A)\label{bound_xdg}
\end{eqnarray}
from
\begin{align}
\E(F)\le a F(A)+1\cdot F(B)=1-(1-a)F(A) \com\qquad 
\E(F)\le \int_A x\rd F(x)+1\cdot F(B)\com\n
\end{align}
respectively.
Next, $\dmint(F,\mu)$ can be written as
\begin{align}
\!\!\!\dmint(F,\mu)
&=
\max_{0\le\nu\le \mom}\E_F[\log(1-(X-\mu)\nu)]\quad \nn
&=
\max_{0\le\nu\le \mom}\left\{
\int_A \log\frac{1-(x-\mu)\nu}{1-(a-\mu)\nu}\rd F(x)+\int_B \log(1-(x-\mu)\nu)\rd F_{(a)}(x)
\right\}.\label{conti_bunkatu}
\end{align}
Since $(1-(x-\mu)\nu)/(1-(a-\mu)\nu)$ is increasing in $\nu$ for $x\le a$,
substituting $0$ and $\momss$ into $\nu$,
we can bound the first term as
\begin{align}
0\,\le\, \int_A \log\frac{1-(x-\mu)\nu}{1-(a-\mu)\nu}\rd F(x)
&\le\,
\int_A \log\frac{1-x}{1-a}\rd F(x)\nn
&\le\,
F(A)\int_A \log(1-x)\frac{\rd F(x)}{F(A)}\quad\;\;\qquad(\mbox{by $a\le 0$})\nn
&\le\,
F(A)\log \left(\int_A(1-x)\frac{\rd F(x)}{F(A)}\right) \qquad(\mbox{Jensen's inequality})\nn
&\le\,
F(A)\log \frac{1-\mt}{F(A)}\per\qquad(\mbox{by \eqref{bound_xdg}})\n
\end{align}
From $\lim_{x\to0}x\log x=0$ and \eqref{bound_ga},
the first term of \eqref{conti_bunkatu} converges to $0$ as $a\to-\infty$.
The second term of \eqref{conti_bunkatu}
equals $\dmint(F_{(a)},\mu)$ and the proof is completed. \fqed
\end{proof1}
Now we show Theorem \ref{thm_bound} based on the preceding lemmas.

\begin{proof2}{of Theorem \ref{thm_bound}}
(i)
The proof is straightforward
since $D(F\Vert G)\ge D(F\Vert G_{(a)})$ always holds for $F\in\cala_a$.

(ii)
First we consider the case that $F$ has a bounded support,
i.e.~$F\in \cala_{a}$ for some $a\in (-\infty,1)$.
It is easily checked that $\dmint(F,\mu)$ defined in \eqref{dmint}
is invariant under the scale transformation
$[0,1]\to [a,1]:x\mapsto a+(1-a)x$. 
Further, $\dmin(F,\mu;\cala_a)$ defined in \eqref{def_dmin} is also invariant
with respect to scale
from the invariance of the divergence.
Since $\dmin(F,\mu;\cala_a)=\dmint(F,\mu)$ 
holds for $a=0$ from Prop.\,\ref{prop_dmin},
it also holds for all finite $a<1$.

Next we consider the case that the support $F$ is not bounded from below.
We show $\dmin\allowbreak (F,\mu;\cala)\le \dmint(F,\mu)$
and $\dmin(F,\mu;\cala)\ge \dmint(F,\mu)$
separately.
We omit the proof for the former part for lack of space,
but it can be proved in a similar procedure as the proof of \citet[Theorem 8]{honda_colt}.

Now we consider the latter inequality.
Take arbitrary $\ep>0$ and
let $a<\mu$ be sufficiently small.
Partitioning $(-\infty,1]$ into $A= (-\infty,a)$
and $B= [a,1]$ we can bound $\dmin(F,\mu;\cala)$ as
\begin{eqnarray}
\inf_{G\in \cala: \E(G)>\mu}D(F\Vert G)
&\ge&
\inf_{G\in \cala: \E(G)>\mu}D(F_{(a)}\Vert G_{(a)})\nn
&\ge&
\inf_{G_{(a)}\in \cala_a: \E(G_{(a)})>\mu}D(F_{(a)}\Vert G_{(a)})\quad(\mbox{by $\E(G)\le\E(G_{(a)})$})\nn
&=&
\dmint(F_{(a)}, \mu)\nn
&\ge&
\dmint(F,\mu)-\ep \quad\mbox{(by Lemma \ref{lem_fa})}\n
\end{eqnarray}
and we complete the proof by letting $\ep\downarrow 0$. \fqed
\end{proof2}
Finally we consider the continuity of $\dmin(F,\mu;\cala)$ in $F$.
\begin{lemma}\label{lem_conti}
If $a<1$ is finite then $\dmin(F,\mu;\cala_a)$ is continuous in $F\in \cala_a$.
\end{lemma}
This lemma is proved for the case $a=0$ in
\citet[Theorem 7]{honda_colt}.
The extension for general bounded supports is straightforward from the scale transformation.

For the case of semi-bounded support distributions, the continuity does not hold any more.
However, we can show the continuity over distributions with expectations
bounded from below.
Here recall that
in view of Theorem \ref{thm_bound}
we write $\dmin(F,\mu)$ instead of
$\dmin(F,\mu;\allowbreak \cala_{-\infty})=\dmint(F,\mu)$
when no confusion arises.
\begin{lemma}\label{lem_levy}
Let $\ep>0$ and $\mu,\mt<1$ be arbitrary.
There exists $\de>0$ such that
\begin{eqnarray}
|\dmin(G,\mu)-\dmin(F,\mu)|\le \ep\label{levy_eq}
\end{eqnarray}
for all $G\in\cala$ such that $\E(G)\ge \mt$ and $\levy(F,G)\le \de$.
\end{lemma}
\begin{proof1}
Applying Lemma \ref{lem_fa} twice to $F$ and $G$,
there exists $a(\ep)$ such that
\begin{eqnarray}
|\dmin(G,\mu)-\dmin(F,\mu)|\le
|\dmin(G_{(a)},\mu)-\dmin(F_{(a)},\mu)|+\ep/2 \label{thm_conti1}
\end{eqnarray}
for all $a\le a(\ep)$ and $G$ such that $\E(G)\ge \mt$.
From the continuity of $\dmin(\cdot,\mu)$ for bounded distribution in Lemma \ref{lem_conti},
there exists $\de(\ep,F_{(a)})$
such that
\begin{eqnarray}
|\dmin(G_{(a)},\mu)-\dmin(F_{(a)},\mu)|\le \ep/2 \label{thm_conti2}
\end{eqnarray}
for all $G_{(a)}$ such that $\levy(G_{(a)}, F_{(a)})\le \de(\ep,F_{(a)})$.
Note that $\levy(G_{(a)},F_{(a)})\le \levy(G,F)$ obviously holds from the definition of
L\'evy distance.
Therefore, from \eqref{thm_conti1} and \eqref{thm_conti2},
we obtain \eqref{levy_eq}
for all $G\in\cala$ such that $\E(G)\ge \mt$
and $\levy(F,G)\le \de(\ep,F_{(a(\ep))})$.\fqed
\end{proof1}%
\section{Large Deviation Probabilities for $\dmin$}
\label{section_LDP}
In this section we consider the behavior of $\dmin(\fhat_t,\mu)$
where $\fhat_t$ is the empirical distribution of $t$ samples
from distribution $F$, which approaches $\dmin(F,\mu)$ as $t$ increases.
For our case of semi-bounded support,
it is sometimes convenient to consider the joint distribution of
empirical mean $\muhat_t=\E(\fhat_t)$ and distribution $\fhat_t$,
since the convergence of the empirical distribution does not mean that of the empirical mean.

Note that, in this section and Appendix \ref{append_ldp},
we sometimes consider moment generating functions
and their Fenchel-Legendre transforms of random variables on domains
other than $\bbR$.
Since the underlying distribution is obvious from the context,
we write e.g.~$\La_{\bbR^2}^*$ to clarify the domain,
whereas the subscript was used to indicate the arm such as $\La_k^*$
in previous sections.
\begin{theorem}\label{ldp_upper}
If $\mu< \E(F)$ and $u>\dmin(F,\mu)$ then
\begin{eqnarray}
P_F[\dmin(\fhat_t,\mu)\ge u \scap \muhat_t\le \mut]
&\le&
\begin{cases}
2\expn{t \La_{\bbR}^*(\mut)}&u\le \La_{\bbR}^*(\mut),\\
2\e(1+t)\expn{tu}&\mathrm{otherwise},
\end{cases}\n
\end{eqnarray}
where $\La_{\bbR}^*(x)\equiv \sup_{\la\in\bbR}\{\la x-\log \E_F[\expp{\la X}] \}$.
\end{theorem}

\begin{theorem}\label{ldp_lower}
Fix arbitrary $\mu>\E(F)$ and $v>0$.
Then it holds for $\cc\ge 2.163$ that 
\begin{eqnarray}
P_F[\dmin(\fhat_t,\mu)\le \dmin(F,\mu)-v]\le \expn{t\uI(v,\,\E(F),\,\mu)}
\n\label{lower_general}
\end{eqnarray}
where 
\begin{eqnarray}
\uI(v,\E(F),\mu)\equiv
\begin{cases}
\frac{v^2}{2(\cc+\frac{1-\E(F)}{1-\mu})}
&v\le \frac12(\cc+\frac{1-\E(F)}{1-\mu})\com\\
\frac{v}{2}-\frac{1}{8}(\cc+\frac{1-\E(F)}{1-\mu})\com
&\mathrm{otherwise.}
\end{cases}\label{ui_def}
\end{eqnarray}

\end{theorem}
We prove these theorems using Prop.\,\ref{cramer}, Theorem \ref{thm_sanov}
and Prop.\,\ref{prop_equiv} in Appendix \ref{append_ldp}.
Before proving Theorem \ref{ldp_upper}, we show its asymptotic version in the following.
\begin{lemma}\label{upper_asymptotic}
If $\mu<\E(F)$ and $u>\dmin(F,\mu)$ then
\begin{eqnarray}
\limsup_{t\to\infty}\frac1t \log P_F[\dmin(\fhat_t,\mu)\ge u\scap \muhat_t\le \mut]
\le -\max\{u, \La_{\bbR}^*(\mut)\}\per\n
\end{eqnarray}
\end{lemma}
\begin{proof1}
Define
$C\equiv\{(G,\E(G)):G\in \cala,\,\dmin(G,\mu)\ge u
\scap \E(G)\le \mut\}\subset \cala\times \bbR$
and let $\bar{C}$ be its closure.
First we show that $\dmin(G,\mu)\ge u$ and $v\le \mu$
for all $(G,v)\in \bar{C}$.

From the definition of closure,
there exists a sequence $\{(G_l,\E(G_l))\in C\}_l$
such that $(G_l,\E(G_l))\to (G,v)$,
i.e., $G_l\to G$ and $\E(G_l)\to v$.
Thus $\E(G_l)\ge v-\ep$ holds for all sufficiently large $l$
where $\ep>0$ is arbitrary.
Therefore, from Lemma \ref{lem_levy} we obtain
\begin{eqnarray}
\dmin(G,\mu)
=\lim_{l\to\infty}\dmin(G_l,\mu)
\ge \liminf_{l\to\infty}u=u\per\n
\end{eqnarray}
The inequality $v\le \mu$ is obvious from $\E(G_l)\to v$ and $\E(G_l)\le \mu$.

Now we obtain from Theorem\,\ref{thm_sanov} that
\begin{eqnarray}
\lefteqn{
\limsup_{t\to\infty}\frac1t \log P_F[\dmin(\fhat_t,\mu)\ge u\scap \muhat_t\le \mut]
}\nn
&\le&\limsup_{t\to\infty}\frac1t \log P_F[(\fhat_t,\muhat_t)\in \bar{C}]\nn
&\le&-\inf_{(G,v)\,:\,\dmin(G,\mu)\ge u \scap v\le \mut}\;\,
\sup_{(\phi,\la)\in \Cb(\bbR)\times \bbR}
\left\{\int \phi(x) \rd G(x)+\la v-\log\int \e^{\phi(x)+\la x}\rd F(x)  \right\}\nn
&\le&-\inf_{(G,v)\,:\,\dmin(G,\mu)\ge u \scap v\le \mu}
\max\{\La_{\bbR}^*(v),\,D(G\Vert F)\}\label{siki_note}\\
&\le&
-\inf_{(G,v)\,:\,\dmin(G,\mu)\ge u \scap v\le \mu}
\max\{\La_{\bbR}^*(v),\dmin(G,\mu)\}
\qquad(\mbox{by $\mu<\E(F)$})\nn
&\le&
-\max\{\La_{\bbR}^*(\mut),u\}\com
\qquad\qquad(\mbox{$\La_{\bbR}^*(v)$ is decreasing in $v\le \mu<\E(F)$})\n
\end{eqnarray}
where \eqref{siki_note} follows from
$(\{0\} \times \bbR)\cup (\Cb(\bbR)\times \{0\})\subset \Cb(\bbR)\times \bbR$
and Prop.\,\ref{prop_equiv}. \fqed
\end{proof1}
\begin{proof2}{of Theorem \ref{ldp_upper}}
Let $\de>0$ be arbitrary and
define $\nu_{i}\equiv 1/(2(1-\mu))+i\de$ for
$i=-M_{\de},-M_{\de}+1,\dots,M_{\de}-1,M_{\de}$,
where $M_{\de}\equiv\left\lfloor1/(2(1-\mu)\de)\right\rfloor$.
Further define $\nu_{-M_{\de}-1}\equiv 0$ and $\nu_{M_{\de}+1}\equiv 1/(1-\mu)$.
Then $\{[\nu_i, \nu_{i+1}]\}$ partitions $[0,\momss]$ into intervals with
length not larger than $\de$.
Therefore the event $\{\dmin(F,\mu)\ge u\}$ can be expressed as
\begin{align}
\lefteqn{
\!\!\!\!\{\dmin(\fhat_t,\mu)\ge u\}
=\left\{\exists \nu\in \left[0,\moms\right],\, L(\nu; \fhat_t,\mu)\ge u\right\}
}\nn
&\!\!\!=\!
\bigcup_{i=-M_{\de}-1}^{-1}\!
\left\{\exists \nu\in \left[\nu_{i},\nu_{i+1}\right], L(\nu; \fhat_t,\mu)\ge u\right\}
\cup\!\bigcup_{i=1}^{M_{\de}+1}
\left\{\exists \nu\in \left[\nu_{i-1},\nu_{i}\right], L(\nu; \fhat_t,\mu)\ge u\right\}.
\label{lde_0}
\end{align}
Since $\nu_{i+1}-\nu_{i}\le \de$ and $L(\nu;\fhat_t,\mu)$ is concave in $\nu$,
it holds for $i\le -1$ that
\begin{align}
\left\{\exists \nu\in \left[\nu_{i},\nu_{i+1}\right],\, L(\nu; \fhat_t,\mu)\ge u\right\}
&\subset
\left\{L(\nu_{i+1}; \fhat_t,\nu)-\de \min\{0,L'(\nu_{i+1};\fhat_t,\nu)\}\ge u\right\}\nn
&\subset
\left\{L(\nu_{i+1}; \fhat_t,\nu)-\de \min\{0,L'(\nu_0;\fhat_t,\nu)\}\ge u\right\}.
\label{lde_1}
\end{align}
Similarly it holds for $i\ge 1$ that
\begin{align}
\left\{\exists \nu\in \left[\nu_{i-1},\nu_{i}\right],\, L(\nu; \fhat_t,\mu)\ge u\right\}
&\subset
\left\{L(\nu_{i-1}; \fhat_t,\mu)+\de \max\{0,L'(\nu_0;\fhat_t,\mu)\}\ge u\right\}. \label{lde_2}
\end{align}

Here the derivative is written as
\begin{eqnarray}
L'(\nu; \fhat_t,\mu)=-\E_{\fhat_t}\left[\frac{X-\mu}{1-(X-\mu)\nu}\right]
=\frac{1}{\nu}-\frac{1}{\nu}\E_{\fhat_t}\left[\frac{1}{1-(X-\mu)\nu}\right]\per\n
\end{eqnarray}
Since $1/(1-(x-\mu)\nu)$ is positive and increasing in $x\le 1$,
it is bounded as
\begin{align}
\frac1{\nu}\ge L'(\nu; \fhat_t,\mu)&\ge\frac1{\nu}-\frac1{\nu}\frac1{1-(1-\mu)\nu}=-\frac{1-\mu}{1-(1-\mu)\nu}\per\n
\end{align}
Thus $L'(\nu_0;\fhat_t,\mu)=L'(1/(2(1-\mu));\fhat_t,\mu)$ is bounded as
\begin{eqnarray}
2(1-\mu)\ge L'(\nu_0;\fhat_t,\mu) \ge -2(1-\mu)\per\n
\end{eqnarray}
Combining this inequality with \eqref{lde_0}, \eqref{lde_1} and \eqref{lde_2}
we obtain
\begin{align}
\lefteqn{
\!\!\!\!\!\! P_{F}[\dmin(\fhat_t,\mu)\ge u\scap \muhat_t\le \mut]
}\nn
&\le
\sum_{-M_{\de}-1\le i\le M_{\de}+1,\,i\neq 0}
P_F\left[L(\nu_i; \fhat_t,\mu)\ge u-2(1-\mu)\de
\scap \muhat_t\le \mut\right].\label{prob_sum}
\end{align}

Now regard $Y= (Y^{(1)},Y^{(2)})\equiv (\log (1-(X-\mu)\nu_i), X)$ as a random variable
on $\mathbb{R}^2$.
Define a closed set $C\equiv [u-2(1-\mu)\de,\,\infty)\times (-\infty,\,\mut]\subset \bbR^2$
and its $\alpha$-blowup
$C^{\alpha}\equiv \allowbreak(u- 2(1-\mu)\de-\alpha,\infty)\times (-\infty,\mut+\alpha)$
for $\alpha>0$.
Then the event $\{L(\nu_i; \fhat_t,\mu)\ge u-2(1-\mu)\de\scap \allowbreak \muhat_t\le \mut\}$
is equivalent to the event that the empirical mean of $Y$ is contained in the
closed convex set $C$.
Thus we obtain from Prop.\,\ref{cramer} (i) that
\begin{eqnarray}
P_{F}[L(\nu_i; \fhat_t,\mu)\ge u-2(1-\mu)\de\scap \muhat_t\le \mut]
&\le&
\exp\left(-t \inf_{y\in C}\La_{\bbR^2}^*(y)\right)\com\label{ldp_upper1}
\end{eqnarray}
where $\La_{\bbR^2}^*(y)$ is defined by \eqref{legendreR}.
Since $C^{\alpha}\supset C$ is open, the exponential rate is bounded as
\begin{eqnarray}
-\inf_{y\in C}\La_{\bbR^2}^*(y)
&\le&
-\inf_{y\in C^{\alpha}}\La_{\bbR^2}^*(y)\nn
&\le&
\liminf_{t\to\infty}\frac1t \log  P_{F}[L(\nu_i; \fhat_t,\mu)> u-2(1-\mu)\de-\alpha
\scap \muhat_t< \mut+\alpha]\nn
&&\phantom{wwwwwwwwwwwwwwwwwwwwwwwwwa}(\mbox{by Prop.\,\ref{cramer} (ii)})\nn
&\le&
\limsup_{t\to\infty} \frac1t\log P_{F}[\dmin(\fhat_t,\mu)\ge u-2(1-\mu)\de-\alpha
\scap \muhat_t\le \mut+\alpha]\nn
&\le&
-\max\{u-2(1-\mu)\de-\alpha,\La_{\bbR}^*(\mut+\alpha)\}\per
\qquad(\mbox{by Lemma \ref{upper_asymptotic}})\n
\end{eqnarray}
Letting $\alpha \downarrow 0$ we obtain
\begin{eqnarray}
-\inf_{y\in C}\La_{\bbR^2}^*(y)\le
-\max\{u-2(1-\mu)\de,\La_{\bbR}^*(\mut)\}
\per \label{ldp_upper2}
\end{eqnarray}

Finally we obtain from \eqref{prob_sum}, \eqref{ldp_upper1} and \eqref{ldp_upper2} that
\begin{eqnarray}
P_F[\dmin(\fhat_t,\mu)\ge x \scap \muhat_t\le \mut]
&\le&
2\left(1+\frac{1}{2(1-\mu)\de}\right)\exp\left(-t\max\{u-2(1-\mu)\de,\La_{\bbR}^*(\mut)\}\right)\n
\end{eqnarray}
and we complete the proof
by letting $\de\to \infty$ for $u\le \La_{\bbR}^*(\mut)$ and
$\de=1/(2t(1-\mu))$ for
$u>\La_{\bbR}^*(\mut)$.
\fqed
\end{proof2}%
\begin{proof2}{of Theorem \ref{ldp_lower}}
Let $u\equiv \dmin(F,\mu)-v$.
First we obtain for $\nu^*=\nu^*(F,\mu)$ that
\begin{eqnarray}
P_F[\dmin(\fhat_t,\mu)\le \dmin(F,\mu)-v]&=&
P_F\left[\max_{0\le \nu\le \momss}\E_{\fhat_t}[\log(1-(X-\mu)\nu)]\le u\right]\nn
&\le&
P_F\big[\E_{\fhat_t}[\log(1-(X-\mu)\nus)]\le u\big]\per\n
\end{eqnarray}
Define random variables $Y\equiv 1-(X-\mu)\nus$ and $Z\equiv \log Y=\log(1-(X-\mu)\nus)$
where $X$ follows the distribution $F$.
Let $\bar{Z}_t$ be the mean of $t$ i.i.d.~copies of $Z$.
Then, from Prop.\,\ref{cramer} (iii), the above probability is bounded as
\begin{eqnarray}
P_F[\dmin(\fhat_t,\mu)\le \dmin(F,\mu)-v]
\,\le\,
P_F[ \bar{Z}_t \le u]
\,\le\,
\expn{t \La_{\bbR}^*(u)}\com\label{owari}
\end{eqnarray}
where
$\La_{\bbR}^*(u)
=\sup_{\la}\big\{\la u-\log \E_F [\e^{\la Z}]\big\}
=\sup_{\la}\big\{\la u-\log \E_F [Y^{\la}]\big\}\n
$. 

Note that $\E_F[\expp{-1\cdot Z}]=\E_F[(1-(X-\mu)\nu^*)^{-1}]\le 1$
from Prop.\,\ref{lem_integrable}
and
$\E_F[\expp{1\cdot Z}]=\E_F[1-(X-\mu)\nus]=1-(\E(F)-\mu)\nus$.
Since they are finite,
the moment generating function $\E_F[\expp{\la Z}]=\E_F[Y^{\la}]$
exists for all $\la\in [-1,1]$
and infinitely differentiable in $\la\in(-1,1)$.

Before evaluating $\La^*(u)$ we bound $\E_F[Y^{\la}]$ for $\la\in[-1,1]$.
%
%
For $\la\in[-1,0]$, we obtain from $\E_F[Y^{-1}]\le 1$
and the convexity of $y^{\la}$ in $\la$ that
\begin{align}
\E_F[Y^{\la}]
\le
\E_F[(-\la)Y^{-1}+(1+\la)Y^0]
\le
-\la+(1+\la)=1\per\label{yla1}
\end{align}
Similarly, we obtain for $\la\in(0,1]$ that
\begin{eqnarray}
\E_F[Y^{\la}]
&\le&
\E_F[(1-\la)Y^0+\la Y^1]\nn
&=&
(1-\la) +\la(1-(\E(F)-\mu)\nus)\nn
&=&
1+\la(\mu-\E(F))\nus\nn
&\le&
1+1\cdot \frac{\mu-\E(F)}{1-\mu}=\frac{1-\E(F)}{1-\mu}\per
\quad\since{by $\mu>\E(F)$ and $\nus\le \mom$}\label{yla2}
\end{eqnarray}
%
%

Define the objective function in $\La_{\bbR}^*(u)$ as
$R(\la)\equiv\la u-\log \E_F[Y^{\la}]$.
Then, for $\la\in[-1,0]$,
\begin{align}
R'(\la)&=u-\frac{\E_F[Y^{\la}\log Y]}{\E_F[Y^{\la}]}
\,\le\, u-\E_F[Y^{\la}\log Y]\per\quad\since{by \eqref{yla1}}\label{lla1}
\end{align}
We bound $R(\la)$ from below for $\la\in [-1/2,0]$ in the following.
For the second term of the right-hand side of \eqref{lla1},
it holds for
$\la\in[-1/2,0]$ that
\begin{eqnarray}
\E_F[Y^{\la}\log Y]&\ge&\E_F[Y^0\log Y]
-\int_{\la}^{0} \max_{\la\in[-\frac12,0]}\left\{\frac{\rd \E_F[Y^{\la}\log Y]}{\rd \la}\right\}\rd \la\nn
&=&
\dmin(F,\mu)+\la \max_{\la \in [-\frac{1}{2},0]}\E_F[Y^{\la}(\log Y)^2]\per\label{mla1}
\end{eqnarray}

Note that $(\log y)^2$ is smaller than $y^{-1/2}$ for $y\to +0$
and smaller than $y$ for $y\to \infty$.
Therefore there exists $\cc>0$ such that
$(\log y)^2\le \cc y^{-1/2} +y$ for all $y>0$.
In fact, this inequality holds by letting $\cc\ge 2.163$.
Then we obtain from \eqref{yla1} and \eqref{yla2} that
\begin{align}
\E_F[Y^{\la}(\log Y)^2]
\le\E_F[Y^{\la}(\cc Y^{-1/2}+Y)]
\le
\cc+\cd\per \label{mla3}
\end{align}
Combining \eqref{lla1}, \eqref{mla1} and \eqref{mla3} with $R(\la)=0$ we obtain
\begin{eqnarray}
R'(\la)&\le& u-\dmin(F,\mu)-\la\left(\cc+\cd\right)
\,=\,-v-\la\left(\cc+\cd\right)\com\nn
R(\la)&=&0+\int_0^{\la} R'(\la)\rd \la
\,\ge\,-\la v-\frac{\la^2}{2}\left(\cc+\cd\right)\per\n
\end{eqnarray}
Finally, 
\begin{align}
\La_{\bbR}^*(u)
\,=\,
\sup_{\la}R(\la)
\,\ge\,
\sup_{\la \in [-\frac{1}{2},0]}R(\la)
\,\ge\,
\begin{cases}
\frac{v^2}{2(\cc+\frac{1-\E(F)}{1-\mu})}\com
&v\le \frac12(\cc+\frac{1-\E(F)}{1-\mu})\com\\
\frac{v}{2}-\frac{1}{8}(\cc+\frac{1-\E(F)}{1-\mu})\com
&\mathrm{otherwise,}
\end{cases}\n
\end{align}
and we obtain the theorem with \eqref{owari}.
\fqed
\end{proof2}

\section{Concluding Remarks}\label{section_remarks}
We proved that the theoretical bound only depends on the upper bound of the
support in the nonparametric stochastic bandits.
We refined the analysis of DMED policy to a non-asymptotic form
for all distributions with moment generating functions in this model.

\appendix
\section{
Large Deviation Principle and its Application to a Joint Distribution
}\label{append_ldp}
\defcitealias{LDP}{DZ}%
In this appendix we consider large deviation principle (LDP) for
the empirical mean and the distribution
based on \cite{LDP} (\citetalias{LDP}, hereafter).
We first summarize results on LDP for empirical means
of finite dimentional random variables and then
we derive LDP for joint distribution of the empirical distribution and the mean
in Theorem \ref{thm_sanov}.

Let $\hat{S}_t$ be the empirical mean of i.i.d.\,random variables
$X_1,\cdots,X_t \in \calX$
with distribution $F$, where $\calX$ is a general topological vector space.
For a distribution on $\bbR$,
we can regard its empirical distribution as the empirical mean
of delta measures $\de_{X_i}\in\cala\subset \calV$,
where $\calV$ is the space of all finite measures on $(-\infty,1]$.
We write
$\muhat_t$ and $\fhat_t$ instead of $\hat{S}_t$ for empirical means
of $X_i\in\bbR$ and $\de_{X_i}\in\cala$, respectively.

Define the logarithmic moment generating function and its Fenchel-Legendre transform
for distribution $F$ by
\begin{eqnarray}
\La_{\calX}(\la)&=&\log \int_{\calX}\e^{\langle \la,u\rangle}\rd F(u)\com\nn
\La_{\calX}^*(x)&=&\sup_{\la\in\calX^*}\{\langle \la,x\rangle-\La(\la)\}\com
\n
\end{eqnarray}
where $\calX^*$ is the space of linear continuous functions on $\calX$.
Especially, for the case $\calX=\bbR^d$
it is expressed for $\calX^*=\bbR^d$ as
\begin{eqnarray}
\langle \la,x\rangle=\sum_{i}\la_i x_i\com\qquad \la,x\in\bbR^d\per\label{legendreR}
\end{eqnarray}
Similarly, for the case $\calX=\calV$, it is expressed for $\calX^*=\Cb(\bbR)$ as
\begin{eqnarray}
\langle \phi,G\rangle=\int \phi(u)\rd G(u)\com \qquad \phi\in \Cb(\bbR),\,G\in\cala\com\n
\end{eqnarray}
where $\Cb(\bbR)$ is the space of
bounded continuous functions on $\bbR$.
Note that it is shown in \citetalias{LDP} that 
in the scope of our paper
$\La_{\calX}^*(x)$ is always a {\it rate function},
that is, a lower semicontinuous function with range $[0,\infty]$,
although we omit this statement in the following.
\begin{proposition}[{\citetalias[Ex.\,2.2.38,\,Theorem 2.2.30 and Lemma 2.2.5]{LDP}}]\label{cramer}
Let $\calX=\bbR^d$
and assume that $\La_{\bbR^d}(\la)$ exists around $\la=0$. 
{\rm (i)} For any convex closed $C \subset\mathbb{R}^{d}$
\begin{eqnarray}
\frac{1}{t}\log P_F[\hat{S}_t\in C]\le-\inf_{x\in C}\La_{\bbR^d}^*(x)\per\n
\end{eqnarray}
{\rm (ii)} For any open $A \subset\mathbb{R}^{d}$
\begin{eqnarray}
\liminf_{t\to\infty}\frac{1}{t}\log P_F[\hat{S}_t\in A]
\ge-\inf_{x\in A}\La_{\bbR^d}^*(x)\per\n
\end{eqnarray}
{\rm (iii)} For the case $d=1$, $\La_{\bbR}^*(x)$ is
decreasing at $x<\E(F)$ and increasing at $x>\E(F)$.
Consequently, 
\begin{eqnarray}
\frac{1}{t}\log P_F[\muhat_t\le x]\le-\La_{\bbR}^*(x)\com &&\quad \mbox{if $x<\E(F)\com$}\nn
\frac{1}{t}\log P_F[\muhat_t\ge x]\le-\La_{\bbR}^*(x)\com &&\quad \mbox{if $x>\E(F)\per$}\n
\end{eqnarray}
\end{proposition}

In well-known Sanov's theorem, LDP
for the empirical distribution is considered.
On the other hand,
in the proof of theorem \ref{ldp_upper},
we have to consider the joint probability that
the empirical distribution and the mean
deviate from a subset of $\cala\times \bbR$.
Theorem \ref{thm_sanov} below is an extension of Sanov's theorem
for this purpose.
This theorem is derived from Cram\'er's theorem
in the same way as the derivation of Sanov's theorem.

Recall that we assume that $\bbR$ is equipped with the standard topology
and $\cala$ is equipped with the topology induced by
L\'evy metric $\levy(F,G)$ for $F,G\in\cala$.
For the space $\cala\times\bbR$ we use the product topology of $\cala$ and $\bbR$,
which is equivalent to the topology induced by the metric
$\max\{\levy(F,G),|x-y|\}$ for $(F,x),\,(G,y)\in \cala\times \bbR$.
\begin{theorem}\label{thm_sanov}
Let $F$ be arbitrary distribution on $\bbR$ such that the moment generating function
exists in some neighborhood of $\la=0$.
For any closed set $C\subset \cala \times \bbR$,
it holds that
\begin{eqnarray}
\limsup_{t\to\infty}\frac1t \log P_F[(\fhat_t,\muhat_t)\in C]\le -\inf_{(G,x)\in C}
\La_{\calV\times \bbR}^*((G,x))\com\label{sanov_eq}
\end{eqnarray}
where 
\begin{eqnarray}
\La_{\calV\times \bbR}^*((G,x))=\sup_{(\phi,\la)\in \Cb(\bbR)\times\bbR}
\left\{\int \phi(u) \rd G(u)+\la x-\log\int \expp{\phi(u)+\la u}\rd F(u)  \right\}\per\n
\end{eqnarray}
\end{theorem}
For the actual computation of $\La_{\calV\times \bbR}^*(\cdot)$
the following proposition is useful.
\begin{proposition}[{\citetalias[Lemma 6.2.13]{LDP}}]\label{prop_equiv}
For all $F,G\in \cala$, 
\begin{eqnarray}
\sup_{\phi\in \Cb(\bbR)}
\left\{\int \phi(u) \rd G(u)-\log\int \expp{\phi(u)}\rd F(u)  \right\}=D(G\Vert F)\per\n
\end{eqnarray}
\end{proposition}

For the rest of this section we prove Theorem \ref{thm_sanov}.
We start with Cram\'er's theorem for general Hausdorff topological
vector spaces $\calX$ and probability measures $F$ on $\calX$.
\begin{proposition}[{\citetalias[Theorem 6.1.3]{LDP}}]\label{cramer_general}
Assume that following (a), (b) hold.
(a) $\calX$  is locally convex
and there exists a closed convex subset $\calE$ of $\calX$ such that $P_F(\calE)=1$.
Further, $\calE$ can be made into a Polish space with respect to
the topology induced by $\calE$.
(b) The closed convex hull of each compact $K\subset \calE$ is compact.
Then it holds for all compact closed set that
\begin{eqnarray}
\limsup_{t\to\infty} \frac1t \log P_F[\hat{S}_t \in C]\le -\inf_{x\in C}\La_{\calX}^*(x)\per
\label{cramer_weak}
\end{eqnarray}
\end{proposition}
The assertion of this proposition is restricted to compact sets
and is called {\it weak} LDP.
We can remove this restriction to {\it full} LDP
if the {\it exponential tightness} is satisfied.
The laws of $\hat{S}_t$ are exponentially tight if,
for every $\alpha<\infty$,
there exists a compact set $K_{\alpha}\subset \calX$ such that
\begin{eqnarray}
\limsup_{t\to\infty}\frac1t\log P_F[\hat{S}_t\in K_{\alpha}^c]<-\alpha\com\n
\end{eqnarray}
where superscript ``{\it c}'' denotes the complement
of the set.
\begin{proposition}[{\citetalias[Lemma 1.2.18]{LDP}}]\label{weak_full}
If the laws of $\hat{S}_t$ are exponentially tight then
\eqref{cramer_weak} holds for all closed set $C$.
\end{proposition}
\begin{proposition}[{\citetalias[Lemma 6.2.6 and Discussion after Eq.\,(2.2.33)]{LDP}}]
\label{tightness}
\noindent
{\rm (i)} The laws of the empirical distributions $\fhat_t\in \cala$ are
exponentially tight for all $F\in \cala$.
{\rm (ii)} The laws of the empirical means $\muhat_t\in \bbR$ are exponentially tight
if the moment generating function $\E_F[\expp{\la X}]$ exists in some neighborhood of
$\la=0$.
\end{proposition}
\begin{proof2}{of Theorem \ref{thm_sanov}}
First we can obtain \eqref{sanov_eq} for all closed compact $C\subset \cala\times \bbR$
as a direct application of Prop.\,\ref{cramer_general} with $\calX:=\calV\times \bbR$
with $\calE:=\cala \times \bbR$ by the following argument.

For the case $\calX:=\calV$ and $\calE:=\cala$,
it is shown as Sanov's theorem that
the assumption of Prop.\,\ref{cramer_general} is satisfied
when $\cala$ is equipped with the topology induced by L\'evy metric
(see \citetalias[Sect.\,6.1]{LDP}).
The essential point in the proof of Sanov's theorem
is that the local convexity in the assumption
is satisfied if a vector space $\calX$ is equipped with
a topology called {\it weak topology}
(see, e.g., \citet[Chap.\,V]{dunford} for detail of weak topologies).
Since the relative topology on $\cala$ of the weak topology of $\calV$
is equivalent to the topology induced by the L\'evy metric,
Prop.\,\ref{cramer_general} is applicable for the case
of Sanov's theorem.
Here note that the weak topology of $\calV\times \bbR$ is equivalent to
the product topology of the weak topologies of $\calV$ and $\bbR$.
Thus it is shown in a parallel way that the assumption is also satisfied
in our case.

In view of Prop.\,\ref{weak_full}, 
we complete the proof if the exponential tightness of the laws of
$(\fhat_t,\,\muhat_t)$ is proved.
From Prop.\,\ref{tightness}, for every $\alpha<\infty$
there exist compact $A_{\alpha}\subset \cala$ and $B_{\alpha}\subset \bbR$
such that
\begin{eqnarray}
\limsup_{t\to\infty}\frac1t\log P_F[\fhat_t\in A_{\alpha}^c]<-\alpha\com \qquad
\limsup_{t\to\infty}\frac1t\log P_F[\muhat_t\in B_{\alpha}^c]<-\alpha\per \label{tight_each}
\end{eqnarray}
Letting $K_{\alpha}:=A_{\alpha}\times B_{\alpha}$ we obtain
\begin{eqnarray}
P_F[(\fhat_t,\,\muhat_t)\in K_{\alpha}^c]\le
P_F[\fhat_t \in A_{\alpha}^c]+P_F[\muhat_t\in B_{\alpha}^c]\per\n
\end{eqnarray}
Combining this inequality
with \eqref{tight_each} we see that the laws of $(\fhat_t,\,\muhat_t)$ are exponentially tight.
\fqed
\end{proof2}

\section{Proof of Theorem \ref{optMD}}\label{proof_opt}
%
%
%
Define events $A_n,\,B_n,\,C_n,\,D_n$ for any $\de>0$ as
\begin{eqnarray*}
A_n&\eqdef&\{\muhat^*(n)\ge \mu^*-\de\}\nn
B_n&\eqdef&
\{\muhat^*(n)\le\mu'+\de\}=
\bigcap_{k=1}^K \{\hat{\mu}_k(n)\le\mu'+\de\}\\
C_n&=&\bigcup_{k\notin \iopt} \{\muhat^*(n)=\muhat_k(n)\ge \mu'+\de\}\\
D_n&=&\bigcup_{k\in \iopt} \{\muhat^*(n)=\muhat_k(n)\le \mu^*-\de\}\per
\end{eqnarray*}
It is easily checked that $\{A_n\cup B_n\cup C_n\cup D_n\}$ is the whole sample space.
Let $\j{n}{i}$ denote the event that arm $i$ is pulled at the $n$-th round
and recall that $\jp{n}{i}$ is given in Algorithm \ref{alg_dmed}.
Then, except for the first $2K$ rounds,
the event $\j{n}{i}$ implies that $\jp{n'}{i}$ occurred for some $K+1\le n'<n$.
Therefore $T_i(n)$ is bounded as
\begin{align}
\lefteqn{T_i(n)}
\nn
&\!=
2+\sum_{t=2}^{\infty} \idxx{\bigcup_{m=2K}^{n-1}
\left\{T_i(m)=t\scap \j{m+1}{i}\right\}}\nn
&\!\le
2+\sum_{t=2}^{\infty}
\idxx{\bigcup_{m=K+1}^{n-1}\!\!\left\{T_i(m)=t\scap \jp{m}{i}\right\}}\nn
&\!\le
2+\sum_{t=2}^{\infty} \idxx{\bigcup_{m=K+1}^{n-1}\!\!\! \left\{T_i(m)=t\scap \jp{m}{i}\cap A_m\right\}}\!
+\sum_{t=2}^{\infty} \idxx{\bigcup_{m=K+1}^{n-1}\!\!\! \left\{T_i(m)=t\scap \jp{m}{i}\cap A_m^c\right\}}\n
\end{align}
and we obtain for the last term that
\begin{align}
\sum_{t=2}^{\infty} \idxx{\bigcup_{m=K+1}^{n-1}\!\! \left\{T_i(m)=t\scap \jp{m}{i}\scap A_m^c\right\}}
&\le
\sum_{m=K+1}^{n-1}\!\!\idx{A_m^c}\nn
&\le
\sum_{m=K+1}^{n-1}\!\!\idx{B_m}+\!\sum_{m=K+1}^{n-1}\!\! \idx{C_m}
+\!\sum_{m=K+1}^{n-1}\!\!\idx{D_m}\,.\n
\end{align}
In the following Lemmas \ref{lemA}--\ref{lemD}
 we bound the expectations of these summations
and they prove the theorem
with $2-2\e/r <0$.
\qed
\begin{lemma}\label{lemA}
Let $i\notin\iopt$ be arbitrary.
If $\xi_{i,\ep,\de}=\ep\di{i}-\de/(1-\mu^*)>0$ then
\begin{align}
\lefteqn{
\phantom{wwwww}
\E\left[
\sum_{t=2}^{\infty}
\idx{\bigcup_{m=K+1}^{n-1} \left\{T_i(m)=t\scap \jp{m}{i}\scap A_m\right\}}
\right]
}\nn
&\phantom{wwwwwwwwww}
\le
\logn+
\frac{1}{1-\expn{\uI(\xi_{i,\ep,\de}\,,\,\mu_i,\,\mu^*)}}\per\n
\end{align}
\end{lemma}

\begin{lemma}\label{lemB}
If $\de<\mu^*-\mu'$ then
\begin{eqnarray}
\E\left[
\sum_{m=K+1}^{n-1} \id[B_m]
\right]
\le
\min_{k\in\iopt}\left\{\frac{2(1+K)}{1-\expn{\mg_k(\tmu)}}
+\frac{2\e}{r\left(1-\expn{r \mg_k(\tmu)}\right)^2}\right\}-\frac{2\e}{r}\per\n
\end{eqnarray}
\end{lemma}

\begin{lemma}\label{lemC}
\begin{eqnarray}
\E 
\left[\sum_{m=K+1}^{n-1} \id[C_m]\right]\le
K\sum_{k\notin \iopt} \frac{1}{1-\expn{\mg_k(\mu'+\de)}}\per\n
\end{eqnarray}
\end{lemma}

\begin{lemma}\label{lemD}
\begin{eqnarray}
\E 
\left[\sum_{m=K+1}^{n-1} \id[D_m]\right]\le
K\sum_{k\in \iopt} \frac{1}{1-\expn{\mg_k(\mu^*-\de)}}\per\n
\end{eqnarray}
\end{lemma}
\begin{proof2}{of Lemma \ref{lemA}}
In the same way as \citet[Lemma 15]{honda_colt},
we obtain
\begin{align}
 \lefteqn{
\!\!\!\!\!
\sum_{t=2}^{\infty}
\idx{\bigcup_{m=K+1}^{n-1} \left\{\jp{m}{i}\scap T_i(m)=t\scap  A_m\right\}}
}\nn
&\!\!\!\le
\logn\nn
&\;\;+\sum_{t=\logn}^{\infty}
\id\left[
\logn (1-r)
\dmin(\fhatt{i}{t},\mu^*-\de)\le \log n 
\right]\nn
&\!\!\!\le
\logn+
\sum_{t=1}^{\infty}
\id\left[
\dmin(\fhatt{i}{t},\mu)\le (1-\ep)\di{i}
\right].
\label{an-kihon}
\end{align}
Note that it holds from Lemma \ref{lem_bibun} and Theorem \ref{ldp_lower} that
\begin{eqnarray}
\lefteqn{
P_{F_i}\left[\dmin(\fhatt{i}{t},\mu^*-\de)\le (1-\ep)\di{i} \right]
}\nn
&\le&
P_{F_i}\left[\dmin(\fhatt{i}{t},\mu^*)-\frac{\de}{1-\mu^*}\le (1-\ep)\di{i}\right]\nn
&\le&
P_{F_i}\left[\dmin(\fhatt{i}{t},\mu^*)\le \di{i}
-\left(\ep\di{i}-\frac{\de}{1-\mu^*}\right)\right]\nn
&\le&
\expn{t \uI(\xi_{i,\ep,\de}\,,\,\mu_i,\,\mu^*)}\per\label{an-p}
\end{eqnarray}
From \eqref{an-kihon} and \eqref{an-p}, we obtain
\begin{align}
\lefteqn{
\!\!\!
\E 
\left[
\sum_{t=2}^{\infty} \idx{\bigcup_{m=K+1}^{n-1}
\left\{T_i(m)=t\scap \jp{i}{m}\scap A_m\right\}}
\right]
}\nn
&\qquad\le
\logn +
\sum_{t=1}^{\infty}
\expn{t\uI(\xi_{i,\ep,\de}\,,\,\mu_i,\,\mu^*)}\nn
&\qquad = \logn+
\frac{1}{1-\expn{\uI(\xi_{i,\ep,\de}\,,\,\mu_i,\,\mu^*)}}.\dqed\n
\end{align}
\end{proof2}%
\begin{proof2}{of Lemma \ref{lemB}}
First we simply bound $\sum_{m=K+1}^{n-1} \id[B_m]$ by
\begin{eqnarray}
\sum_{m=K+1}^{n-1} \id[B_m]
&\le&
\sum_{t=1}^{\infty} \sum_{m=K+1}^{\infty}
\id[B_m \scap T_k(m)=t]\com\label{bound_bn}
\end{eqnarray}
where $k\in\iopt$ is arbitrary.
By the same argument as \citet[Lemma 16]{honda_colt},
the event
$\{\dmin(\fhatt{k}{t},\tmu)\le u \scap \muhatt{k}{t}\le \tmu\}$ implies
\begin{eqnarray}
\sum_{m=K+1}^{\infty}\id[B_m \scap T_k(m)=t]\le \expp{tu(1-\para)}+K\per\n
\end{eqnarray}

Let $P(u)\equiv P_{F_k}[\dmin(\fhatt{k}{t},\tmu)\ge u \scap \muhatt{k}{t}\le \tmu]$.
When we simply write $\mg_k$ for $\mg_k(\tmu)$ given in \eqref{def_ik},
it holds from Theorem \ref{ldp_upper} that
\begin{eqnarray}
\lefteqn{
\E\left[\sum_{m=K+1}^{\infty}\id[B_m \scap T_k(m)=t]\right]
}\nn
&\le&
\int_{\infty}^{0} (\expp{tu(1-r)}+K)\rd P(u)\nn
&=&
\left[(\expp{tu(1-\para)}+K)P(u)\right]_{\infty}^{0}
+t(1-r)\int_0^{\infty}\expp{tu(1-r)}P(u)\rd u
\nn
&\le&
2(1+K)\expn{t \mg_k}
+2t(1-r)\int_0^{\mg_k}
\expn{t(\mg_k-(1-r)u)}
\rd u
+2\e t(1-r)(1+t)\int_{\mg_k}^{\infty}
\expn{tru}
\rd u\nn
&\le&
2(1+K)\expn{t \mg_k}
+2\expn{tr \mg_k}
+\frac{2\e (1-r)}{r}
(1+t)\expn{tr \mg_k}\nn
&\le&
2(1+K)\expn{t \mg_k}
+\frac{2\e}{r}
(1+t)\expn{tr \mg_k}\per\n
\end{eqnarray}
Taking the summation over $t$ with formula
\begin{align}
\sum_{t=1}^{\infty}(1+t)\rho^t\,=\,\frac{1}{(1-\rho)^2}-1\com\n
\end{align}
we obtain from \eqref{bound_bn} that
\begin{eqnarray}
\E_F\left[\sum_{m=K+1}^{n-1}\id[B_m]\right]
&\le&
\frac{2(1+K)}{1-\expn{\mg_k}}
+\frac{2\e}{r\left(1-\expn{r \mg_k}\right)^2}-\frac{2\e}{r}\per\label{bn_last}
\end{eqnarray}
We complete the proof by taking $k\in\iopt$ such that \eqref{bn_last} is minimized.
\fqed
\end{proof2}%
\begin{proof2}{of Lemmas \ref{lemC} and \ref{lemD}}
We obtain from the definition of $C_n$ that
\begin{eqnarray}
\sum_{m=K+1}^{n-1} \id[C_m]
&\le&
\sum_{k\notin \iopt} \sum_{m=K+1}^{\infty} \id[\muhat^*(m)=\muhatn{k}{m}\ge \mu'+\de
]\nn
&\le&
\sum_{k\notin \iopt} \sum_{t=1}^{\infty}\sum_{m=K+1}^{\infty}
 \id[\muhat^*(m)=\muhatt{k}{t}\ge \mu'+\de \scap T_k(m)=t]\per\label{lemc1}
\end{eqnarray}
By the same argument as \citet[Lemma 17]{honda_colt},
we have
\begin{eqnarray}
\sum_{m=K+1}^{\infty}
 \id[\muhat^*(m)=\muhatt{k}{t} \scap T_k(m)=t]\le K\per\label{lemc2}
\end{eqnarray}
On the other hand, from Prop.\,\ref{cramer} (iii) we have
\begin{eqnarray}
P_{F_k}[\muhatt{k}{t}\ge \mu'+\de]
\le
\expn{t \mg_k(\tmu)}\com \label{lemc3}
\end{eqnarray}
where $\mg_k(x)$ is given in \eqref{def_ik}.
Finally we obtain from \eqref{lemc1}--\eqref{lemc3} that
\begin{eqnarray}
\E\left[\sum_{m=K-1}^{n-1} \id[C_m] \right]
\,\le\,
K\sum_{k\notin \iopt} \sum_{t=1}^{\infty}
 P_{F_k}[\muhatt{k}{t}\ge \mu'+\de]
\,\le\,
K\sum_{k\notin \iopt} \frac{1}{1-\expn{\mg_k(\mu'+\de)}}
\n
\end{eqnarray}
and Lemma \ref{lemC} is proved.
In the same way,
we obtain Lemma \ref{lemD} from
\begin{align}
\E 
\left[\sum_{m=K-1}^{n-1} \id[D_m] \right]
 &\le
K\sum_{k\in \iopt} \sum_{t=1}^{\infty}
 P_{F_i}[\muhatt{k}{t}\le \mu^*-\de]\nn
  &\le
K\sum_{k\notin \iopt} \sum_{t=1}^{\infty}\expn{t\mg_k(\mu^*-\de)}\nn
&\le
K\sum_{k\notin \iopt} \frac{1}{1-\expn{\mg_k(\mu^*-\de)}}\per
\dqed\n
\end{align}
\end{proof2}




\end{document}